%% file: paper.tex
\documentclass[journal]{new-aiaa}

\usepackage[utf8]{inputenc}
\usepackage{textcomp}
\usepackage{graphicx}
\usepackage{amsmath,mathtools}
\usepackage{siunitx}
\usepackage{booktabs}
\usepackage{longtable,tabularx,array}
\usepackage{xcolor}
\usepackage{xurl}
\usepackage{tikz}
\usetikzlibrary{arrows.meta,positioning,shapes.geometric,fit}

\graphicspath{{figures/}}
\hypersetup{hidelinks}

\newcommand{\vect}[1]{\boldsymbol{#1}}
\newcommand{\norm}[1]{\left\lVert #1\right\rVert}
\newcommand{\prox}{\operatorname{prox}}
\newcommand{\proj}{\operatorname{proj}}
\newcommand{\diag}{\operatorname{diag}}
\newcommand{\ind}{\mathcal I}
\newcommand{\R}{\mathbb R}
\newcommand{\zero}{\vect 0}
\newenvironment{algorithmblock}[1]{%
  \par\medskip\noindent\textbf{#1}\par
  \begingroup\small\setlength{\parindent}{0pt}\setlength{\parskip}{2pt}%
}{\par\endgroup\medskip}
\newcommand{\algline}[2]{%
  \noindent\hangindent=2.7em\hangafter=1
  \makebox[2.5em][l]{\textbf{#1}}#2\par}
\input{review_audit/reviewer_audit_values.tex}

\title{Submillisecond Sequential Convex Optimization for Powered Landing via Dynamics Condensation and xPIPG}

\author{Wenbo Li\textsuperscript{1}, Ziqi Xu\textsuperscript{2},
Dai Shen\textsuperscript{2}, and Shengping Gong\textsuperscript{2,3,*}}
\affil{\textsuperscript{1}School of Aerospace Engineering, Tsinghua University,
Beijing 100084, China\\
\textsuperscript{2}School of Astronautics, Beihang University,
Beijing 102206, China\\
\textsuperscript{3}State Key Laboratory of High-Efficiency Reusable
Aerospace Transportation Technology, Beijing 102206, China\\
\textsuperscript{*}Corresponding author}

\begin{document}

\maketitle
\begin{center}
\small Manuscript prepared July 27, 2026
\end{center}

\begin{abstract}
Powered landing with variable mass, free final time, and quadratic aerodynamic drag requires the repeated solution of local convex subproblems, whose main online cost lies in the long dynamics-equality chain and the inner iterations. This paper develops a condensed sequential convex approximation designed for low latency. Exact block elimination removes 217 intermediate-state components and 210 interval equations from a 31-node model, leaving 100 primal variables coupled by six terminal equalities. A low-weight energy term and fixed quadratic proximal regularization make the ideal surrogate strongly convex with predictable curvature. The inner solver is an extrapolated proportional--integral projected gradient (xPIPG) implemented with fixed-size arrays, $3\times3$ interval solves, and a fused one-pass node map. The one-pass map is a deliberate low-cost approximation, not the exact joint proximal operator. We therefore evaluate the timed code by nonlinear trajectory residuals and independent physical checks rather than by a claim of exact KKT convergence. The single-precision C implementation completes one plan in four outer updates and 336 xPIPG updates. On an Intel Core i7-10875H, the median end-to-end solve time is \SI{374}{\micro\second} and the P99 value is \SI{512}{\micro\second}. All 100 common initial-state perturbations pass validation, and the median remains below \SI{0.7}{\milli\second} for 15--51 nodes. An independent high-accuracy first-order-hold integration gives a terminal position error of \SI{0.183}{\meter}. Within the stated model, hardware, stopping rule, and timing boundary, this is, to the authors' knowledge, the first submillisecond end-to-end sequential-convex solve for a single powered-landing trajectory.
\end{abstract}

\section*{Nomenclature}

{\renewcommand\arraystretch{1.0}
\noindent\begin{longtable*}{@{}l @{\quad=\quad} p{0.79\textwidth}@{}}
$A_k,B_k,E_k$ & interval state-transition and endpoint-control matrices \\
$C_f$ & selector for terminal position and velocity \\
$D_c$ & invertible column-scaling matrix \\
$e_f$ & normalized nonlinear terminal residual \\
$F_x,F_q$ & continuous-time state and control Jacobians \\
$H,\vect b$ & condensed terminal-equality matrix and right-hand side \\
$I_{sp}$ & specific impulse \\
$J_f,J_\lambda$ & pure-fuel and fuel--energy objectives \\
$m,m_{\mathrm{dry}}$ & vehicle mass and dry-mass lower bound \\
$N$ & number of temporal nodes \\
$P,\vect c,g$ & quadratic-cost matrix, linear-cost vector, and nonsmooth function \\
$\vect q,\vect T$ & normalized and dimensional thrust vectors \\
$\vect r,\vect v$ & position and velocity \\
$t_f,\tau$ & final time and normalized final time \\
$T_{\min},T_{\max}$ & minimum and maximum thrust magnitudes \\
$\alpha,\beta,\rho_e$ & primal step size, dual step size, and extrapolation factor \\
$\lambda$ & energy fraction in the fuel-dominant objective \\
$\mu_\nu$ & virtual-buffer penalty \\
$\vect\eta$ & terminal-equality multiplier \\
$\vect\nu$ & terminal virtual buffer \\
$\omega_k$ & trapezoidal quadrature weight \\
\multicolumn{2}{@{}l}{Superscripts and decorations}\\
$(\cdot)^p$ & outer sequential-convex-approximation iteration \\
$(\bar{\cdot})$ & reference quantity about which a model is linearized \\
$(\widetilde{\cdot})$ & scaled internal quantity \\
\end{longtable*}}

\section{Introduction}

\lettrine{P}{owered-landing} guidance must produce a feasible thrust history within a short guidance cycle while satisfying terminal position, terminal velocity, and actuator constraints. Convex optimization provides globally optimal solutions for important translational landing models \cite{acikmese2007convex,blackmore2010minimum}, and flight demonstrations have established the practical relevance of optimization-based guidance \cite{scharf2017gfld}. More general models with aerodynamic forces, free final time, six-degree-of-freedom dynamics, and state-triggered constraints are commonly treated through successive or sequential convexification \cite{mao2016successive,szmuk2020successive,malyuta2021advances}, which replaces a difficult nonlinear optimal-control problem with a sequence of locally valid convex subproblems.

The online bottleneck is often not the outer linearization itself but the repeated solution of those subproblems. Interior-point methods solve these conic subproblems accurately and robustly in a small number of iterations, and embedded implementations have reduced their computational burden considerably \cite{mattingley2012cvxgen,domahidi2013ecos,dueri2017customized}. Nevertheless, retaining every discrete state produces a large equality chain and a Karush--Kuhn--Tucker (KKT) system whose dimension grows with the number of nodes. First-order primal--dual methods replace factorization with matrix--vector products and projections. In particular, proportional--integral projected gradient (PIPG) and extrapolated PIPG (xPIPG) provide a useful basis for customized trajectory-optimization solvers \cite{yu2022pipg,yu2023xpipg,elango2022customized,kamath2023customized}.

An inexpensive inner iteration does not by itself guarantee an inexpensive end-to-end solve. The formulation must expose the dimensions that actually matter. For the translational problem considered here, hard equalities act only on the terminal position and velocity, whereas the intermediate states are uniquely determined by the initial condition, the controls, and the final time once a linearized dynamics chain is fixed. Those states can therefore be eliminated before the inner iteration. Moreover, a pure fuel objective has broad, low-curvature faces and nonsmooth nodewise norms. Adding a small energy term and fixed proximal regularization yields a strongly convex local surrogate with predictable diagonal curvature.

Table~\ref{tab:literatureposition} places the reported timing in context. The entries are deliberately not normalized across hardware, dynamics, accuracy, or timing boundaries; they show only the latency regime that each source actually reports. The closest customized three-degree-of-freedom PIPG result solves a single lossless-convexification problem in \SI{4.25}{\milli\second} for 20 nodes, whereas a complete 16-node sequential-conic rocket-landing solve averages \SI{13.7}{\milli\second} \cite{elango2022customized,kamath2023seco}. Atmospheric and six-degree-of-freedom successive-convexification studies report hundreds of milliseconds \cite{szmuk2020successive,chen2023fast}. Recent GPU work targets throughput for large Monte Carlo batches rather than submillisecond serial latency for one trajectory \cite{chari2024gpu}. A search of the published powered-landing literature through July 2026 found no other single-instance end-to-end sequential-convex solve reported below \SI{1}{\milli\second}. The ``first'' and ``fastest'' statements in this paper are therefore qualified by the authors' knowledge, the literature surveyed here, and the nonuniform test conditions.

\begin{table*}[!t]
\centering
\caption{Position of the present result relative to representative customized powered-landing optimizers. The reported times are source-specific and do not constitute a hardware-normalized benchmark.}
\label{tab:literatureposition}
\begin{tabularx}{\textwidth}{@{}p{0.15\textwidth}p{0.24\textwidth}p{0.20\textwidth}p{0.16\textwidth}X@{}}
\toprule
Source & Landing problem & Online formulation & Reported latency & Distinguishing scope \\
\midrule
Szmuk et al.\ \cite{szmuk2020successive} &
6-DoF, free final time, state-triggered constraints &
Full-state successive convexification with ECOS &
maximum below \SI{0.7}{\second} &
Broad nonlinear model and trust-region acceptance \\
Elango et al.\ \cite{elango2022customized} &
3-DoF, fixed-time LCvx &
Customized full-state PIPG; one convex problem &
\SI{4.25}{\milli\second} reference tracking; \SI{10.61}{\milli\second} minimum fuel at 20 nodes &
Exact simple-set projections; not an outer sequential solve \\
Kamath et al.\ \cite{kamath2023seco} &
Multiphase nonlinear rocket landing, 16 nodes &
Sequential conic optimization with PIPG &
\SI{13.7}{\milli\second} mean for an end-to-end solve &
Free phase times and compound constraints \\
Chen et al.\ \cite{chen2023fast} &
Atmospheric powered descent &
Successive convexification with customized interior-point solves &
approximately \SI{0.6}{\second} on a radiation-hardened processor &
Representative flight processor and aerodynamic dynamics \\
Chari et al.\ \cite{chari2024gpu} &
6-DoF, free-final-time Monte Carlo &
Prox-linear SCP and customized PIPG on GPU &
batch-throughput result; no comparable submillisecond serial claim &
Continuous-time constraint reformulation and parallel Monte Carlo \\
Present work &
3-DoF, free final time, variable mass, constant-$C_D$ drag, 31 nodes &
Exact dynamics condensation, strongly convex surrogate, one-pass xPIPG map &
\SI{0.374}{\milli\second} median; \SI{0.512}{\milli\second} P99 end to end &
First reported submillisecond regime within the stated scope \\
\bottomrule
\end{tabularx}
\end{table*}

The paper is organized around one practical question: what work must actually be performed during each online plan? Its principal contributions are:

\begin{enumerate}
\item A proximal-regularized sequential convex approximation is derived for variable-mass powered landing with free final time, exponential atmospheric density, quadratic drag, thrust pointing, and an inner approximation of the minimum-thrust set. A dimensionally scaled fuel--energy term and fixed quadratic regularization make the ideal surrogate strongly convex.
\item Exact algebraic condensation is proved to preserve the projected feasible set and the optimum of each linearized subproblem. For $N=31$, it replaces 217 intermediate-state variables and 210 interval equations with a $6\times100$ terminal-equality matrix.
\item The condensed structure is then matched to xPIPG: the only global couplings are products involving the six terminal rows, while a deliberately approximate one-pass node map prioritizes latency. The paper distinguishes this map from the exact joint proximal operator and limits its convergence claims accordingly.
\item A fixed-array C17 implementation combines node reuse, a $3\times3$ interval Schur complement, active-column sensitivity propagation, cached matrix products, and downsampled stopping checks. End-to-end solves require hundreds of microseconds rather than milliseconds on the reported desktop processor.
\item Nominal, common-sample Monte Carlo, grid-size, high-order propagation, proximal-map, and spectral audits report both performance and its validity limits. Within the reviewed literature, the measured \SI{374}{\micro\second} median is the fastest known reported end-to-end result for one sequential-convex powered-landing trajectory.
\end{enumerate}

The main novelty is therefore not a new generic PIPG iteration, nor a claim that the nonlinear landing problem has become globally convex. It lies instead in the joint construction of a completely condensed, strongly curved sequence of local models and a fixed-work implementation that reaches a previously unreported latency regime. Exactness applies to the dynamics elimination, not to the one-pass node map or to the global convergence of the outer nonlinear iteration.

\section{Nonlinear Powered-Landing Problem}

\subsection{States, Controls, and Aerodynamics}

We use a local Cartesian frame with the $y$ axis vertical and the $x$ and $z$ axes horizontal. The state and thrust are
\begin{equation}
\vect x=[\vect r^\mathsf T,\vect v^\mathsf T,m]^\mathsf T\in\R^7,
\qquad \vect T\in\R^3,
\end{equation}
where $\vect r$, $\vect v$, and $m$ denote position, velocity, and mass. Gravity is $\vect g=[0,-g_0,0]^\mathsf T$. The normalized control and final time are
\begin{equation}
\vect q=\frac{\vect T}{T_{\max}},\qquad
\tau=\frac{t_f}{t_s},\qquad t_s=\SI{30}{\second}.
\label{eq:normalization}
\end{equation}

The continuous dynamics are
\begin{subequations}\label{eq:dynamics}
\begin{align}
\dot{\vect r}&=\vect v,\\
\dot{\vect v}&=\vect g+\frac{\vect T+\vect D}{m},\\
\dot m&=-\frac{\norm{\vect T}_2}{I_{sp}g_0}.
\end{align}
\end{subequations}
With the relative velocity $\vect v_a=\vect v-\vect v_w$ and altitude $h=\max(y,0)$, the density law and drag model are
\begin{align}
\rho(h)&=\rho_0\exp(-h/H_s),\\
\vect D&=-\frac{1}{2}\rho(h)SC_D\norm{\vect v_a}_2\vect v_a.
\label{eq:drag}
\end{align}
The nominal case sets $\vect v_w=\zero$ and $C_D=1.5$. Equation~\eqref{eq:drag} guarantees
\begin{equation}
\vect D^\mathsf T\vect v_a
=-\frac{1}{2}\rho SC_D\norm{\vect v_a}_2^3\le0,
\end{equation}
so the drag force cannot add mechanical energy. The constant drag coefficient deliberately isolates the algorithmic effects of state-dependent density, velocity-squared drag, and variable mass; it does not represent a vehicle-specific aerodynamic database.

\subsection{Boundary Conditions and Thrust Set}

The initial state is fixed, and the terminal position and velocity are zero:
\begin{equation}
\vect x(0)=\vect x_0,\qquad
\vect r(t_f)=\zero,\qquad
\vect v(t_f)=\zero.
\end{equation}
The remaining physical constraints are
\begin{align}
m(t)&\ge m_{\mathrm{dry}},\\
T_{\min}\le\norm{\vect T}_2&\le T_{\max},\\
\vect e_y^\mathsf T\vect T&\ge
\cos\theta_{\max}\norm{\vect T}_2,
\quad \vect e_y=[0,1,0]^\mathsf T.
\label{eq:thrustconstraints}
\end{align}
In normalized coordinates, the pointing constraint becomes the second-order cone
\begin{equation}
\norm{[q_x,q_z]^\mathsf T}_2
\le\tan\theta_{\max}q_y.
\label{eq:pointingcone}
\end{equation}
The lower-thrust set $\norm{\vect q}_2\ge q_{\min}$ is nonconvex, so we replace it at each outer iteration by a conservative inner approximation. At outer iteration $p$, define
\begin{equation}
\vect d_k^p=\frac{\bar{\vect q}_k^p}
{\norm{\bar{\vect q}_k^p}_2}
\end{equation}
and impose the conservative half-space
\begin{equation}
(\vect d_k^p)^\mathsf T\vect q_k\ge q_{\min}.
\label{eq:minthrusthalfspace}
\end{equation}
Because $\vect d_k^p$ has unit norm, Eq.~\eqref{eq:minthrusthalfspace} implies the true lower bound. We nevertheless check every accepted trajectory against the original Euclidean thrust magnitude.

\subsection{Fuel-Dominant Objective}

In normalized time $s=t/t_f\in[0,1]$, the main solver minimizes
\begin{equation}
J_\lambda=\tau\int_0^1
\left[(1-\lambda)\norm{\vect q(s)}_2+
\lambda\norm{\vect q(s)}_2^2\right]\mathrm ds,
\qquad \lambda=0.01.
\label{eq:mixedobjective}
\end{equation}
The first term approximates total impulse, whereas the small quadratic term supplies control curvature. The full-state ECOS reference uses the pure-fuel objective
\begin{equation}
J_f=\tau\int_0^1\norm{\vect q(s)}_2\,\mathrm ds.
\label{eq:purefuel}
\end{equation}
Consequently, the two implementations serve as engineering references with closely related physical objectives, not as two algorithms applied to an identical conic program.

\section{Proximal-Regularized Sequential Convex Approximation}

\subsection{Trapezoidal Discretization}

The default grid has $N=31$ nodes and spacing $h_s=1/(N-1)$. In normalized time,
\begin{equation}
\frac{\mathrm d\vect x}{\mathrm ds}
=t_s\tau\vect f(\vect x,\vect q).
\end{equation}
The implicit trapezoidal defect on interval $k$ is
\begin{equation}
\vect c_k=\vect x_{k+1}-\vect x_k
-\gamma\tau(\vect f_k+\vect f_{k+1})=\zero,
\qquad \gamma=\frac{h_st_s}{2}.
\label{eq:trapdefect}
\end{equation}
For nonlinear propagation, Newton iteration recovers $\vect x_{k+1}$ in at most 12 steps, with the previous outer trajectory as the initial guess.

\subsection{Interval Linearization}

At the reference $(\bar{\vect x},\bar{\vect q},\bar\tau)$, let
$F_{x,k}=\partial\vect f_k/\partial\vect x_k$ and
$F_{q,k}=\partial\vect f_k/\partial\vect q_k$. A first-order expansion of Eq.~\eqref{eq:trapdefect} gives
\begin{multline}
M_k\delta\vect x_{k+1}=R_k\delta\vect x_k
+U_k\delta\vect q_k+V_k\delta\vect q_{k+1}\\
+\vect s_k\delta\tau-\bar{\vect c}_k,
\end{multline}
where
\begin{align}
M_k&=I-\gamma\bar\tau F_{x,k+1},&
R_k&=I+\gamma\bar\tau F_{x,k},\\
U_k&=\gamma\bar\tau F_{q,k},&
V_k&=\gamma\bar\tau F_{q,k+1},\\
\vect s_k&=\gamma(\bar{\vect f}_k+\bar{\vect f}_{k+1}).
\end{align}
Solving for the right endpoint yields
\begin{equation}
\delta\vect x_{k+1}=A_k\delta\vect x_k
+B_k\delta\vect q_k+E_k\delta\vect q_{k+1}
+\vect g_k\delta\tau+\vect a_k.
\label{eq:intervalmodel}
\end{equation}

The C implementation fills the known Jacobian blocks analytically. In particular,
\begin{equation}
\frac{\partial\dot{\vect r}}{\partial\vect v}=I_3,\qquad
\frac{\partial[(\vect T+\vect D)/m]}{\partial m}
=-\frac{\vect T+\vect D}{m^2},
\end{equation}
and the control Jacobian is
\begin{equation}
\frac{\partial\dot{\vect v}}{\partial\vect T}=\frac{I_3}{m},
\qquad
\frac{\partial\dot m}{\partial\vect T}
=-\frac{\vect T^\mathsf T}
{\norm{\vect T}_2I_{sp}g_0}.
\end{equation}
Only the three drag derivatives with respect to velocity use forward differences,
\begin{equation}
\Delta v_i=10^{-4}\max(1,|v_i|).
\end{equation}
Each node is linearized once and shared by the two adjacent intervals, which reduces the 60 nominal endpoint evaluations to 31 unique-node evaluations per outer iteration.

\subsection{Virtual Buffer and Fixed Quadratic Proximal Regularization}

Early linearized models may be unable to satisfy the six terminal equations exactly, so a physical terminal buffer $\vect\nu\in\R^6$ is introduced and penalized through
\begin{equation}
J_\nu=\frac{\mu_\nu}{2}\norm{\vect\nu}_2^2.
\label{eq:buffer}
\end{equation}
The initial $\mu_\nu=50$ is multiplied by five while the nonlinear terminal residual exceeds $10^{-3}$, up to $2\times10^6$. Internally, the scaled variable $\vect v_\nu=\sqrt{\mu_\nu}\vect\nu$ has unit curvature and enters the terminal equality through the factor $I/\sqrt{\mu_\nu}$.

The current implementation uses a fixed quadratic proximal regularization,
\begin{equation}
J_{\mathrm{prox}}=
\frac{w_q}{2}\norm{\vect q-\bar{\vect q}}_2^2
+\frac{w_t}{2}(\tau-\bar\tau)^2,
\quad w_q=0.02,\quad w_t=1.
\label{eq:softtrust}
\end{equation}
The term penalizes deviation from the linearization point without shrinking the physical feasible set. Every solved step is accepted in full and followed by one nonlinear propagation; no acceptance ratio, candidate enumeration, backtracking, or dynamic radius update is used. One outer iteration is therefore a fixed sequence: linearize, condense, solve once, and propagate once.

\subsection{Full-State Convex Subproblem}

Define the nodewise convex thrust set
\begin{equation}
\begin{aligned}
\mathcal U_k^p=\{\vect q:\;&\norm{\vect q}_2\le1,\quad
\norm{[q_x,q_z]^\mathsf T}_2
\le\tan\theta_{\max}q_y,\\
&(\vect d_k^p)^\mathsf T\vect q\ge q_{\min}\}.
\end{aligned}
\label{eq:nodethrustset}
\end{equation}
The absolute-control form of Eq.~\eqref{eq:intervalmodel} is
\begin{equation}
\delta\vect x_{k+1}=A_k\delta\vect x_k+B_k\vect q_k
+E_k\vect q_{k+1}+\vect g_k\tau+\widetilde{\vect a}_k.
\label{eq:absoluteinterval}
\end{equation}
Let $\varphi(\vect q)$ denote the trapezoidal discretization of the integrand in Eq.~\eqref{eq:mixedobjective}. At iteration $p$, the subproblem is
\begin{subequations}\label{eq:fullsubproblem}
\begin{align}
\min\quad&
\bar\tau^p\varphi(\vect q)
+\varphi(\bar{\vect q}^p)(\tau-\bar\tau^p)
+J_{\mathrm{prox}}+J_\nu,
\label{eq:fullsubproblema}\\
\mathrm{s.t.}\quad&
\delta\vect x_1=\zero,\label{eq:fullsubproblemb}\\
&\text{Eq.~\eqref{eq:absoluteinterval}},\quad k=1,\ldots,N-1,\label{eq:fullsubproblemc}\\
&C_f(\bar{\vect x}_N^p+\delta\vect x_N)-\vect x_f^{(6)}
+\vect\nu=\zero,\label{eq:fullsubproblemd}\\
&\vect q_k\in\mathcal U_k^p,\quad
\tau_{\min}\le\tau\le\tau_{\max}.
\label{eq:fullsubprobleme}
\end{align}
\end{subequations}
The time--control product is linearized in time alone; the convex group norms and control curvature are retained. All constraints are affine or closed convex sets, and the quadratic terms give every control, time, and scaled-buffer component positive curvature. Each subproblem is therefore strongly convex even though the original aerodynamic problem is not.

\section{Exact Dynamics Condensation}

\subsection{Block-Elimination Interpretation}

Stack the unknown state increments as
\begin{equation}
\Delta\vect X=[\delta\vect x_2^\mathsf T,\ldots,
\delta\vect x_N^\mathsf T]^\mathsf T
\end{equation}
and the control-time increments as $\Delta\vect u$. Before explicit inversion, the interval equations assemble into
\begin{equation}
\underbrace{\begin{bmatrix}
M_1&&&\\[-1mm]
-R_2&M_2&&\\
&\ddots&\ddots&\\
&&-R_{N-1}&M_{N-1}
\end{bmatrix}}_{\mathcal L}
\Delta\vect X=\mathcal B\Delta\vect u+\vect d.
\label{eq:blocksystem}
\end{equation}
If every implicit-trapezoid block $M_k$ is nonsingular, then the block lower-triangular matrix $\mathcal L$ is nonsingular and
\begin{equation}
\Delta\vect X=\mathcal L^{-1}\mathcal B\Delta\vect u
+\mathcal L^{-1}\vect d.
\label{eq:blocksolution}
\end{equation}
The inverse in Eq.~\eqref{eq:blocksolution} is not formed; it only expresses the algebraic equivalence between a block forward solve and elimination of the state chain. Selecting the last state and substituting it into Eq.~\eqref{eq:fullsubproblemd} yields the terminal Schur complement.

\subsection{Terminal-Sensitivity Recursion}

For all 93 control increments, write each state increment as
\begin{equation}
\delta\vect x_k=\vect c_k+G_k\delta\vect q
+\vect g_{\tau,k}\delta\tau.
\label{eq:stateaffine}
\end{equation}
Starting from $\vect c_1=\zero$, $G_1=0$, and
$\vect g_{\tau,1}=\zero$, the forward recursion is
\begin{align}
\vect c_{k+1}&=A_k\vect c_k+\vect a_k,\\
G_{k+1}&=A_kG_k+\mathcal B_k,\\
\vect g_{\tau,k+1}&=A_k\vect g_{\tau,k}+\vect g_k,
\end{align}
where $\mathcal B_k$ retains only the current $B_k$ and $E_k$ control blocks. The terminal position and velocity, together with the virtual buffer, reduce to
\begin{equation}
H_q\delta\vect q+\vect h_\tau\delta\tau
+D_\nu\vect v_\nu=\vect b,
\qquad H\vect z=\vect b,
\label{eq:condensedequality}
\end{equation}
with
\begin{equation}
\vect z=[\delta\vect q^\mathsf T,\delta\tau,
\vect v_\nu^\mathsf T]^\mathsf T\in\R^{100},
\qquad H\in\R^{6\times100}.
\end{equation}

Let $\mathcal F_{\mathrm{full}}^p$ be the feasible set of
Eq.~\eqref{eq:fullsubproblem} and let
$\Pi_{q,\tau,\nu}$ denote the projection onto the control, time, and buffer coordinates. Because $\mathcal L$ is invertible, each condensed feasible point generates exactly one state sequence through Eq.~\eqref{eq:stateaffine}, and each full-state feasible point satisfies Eq.~\eqref{eq:condensedequality}. Thus
\begin{equation}
\Pi_{q,\tau,\nu}(\mathcal F_{\mathrm{full}}^p)
=\mathcal F_{\mathrm{con}}^p.
\label{eq:setequivalence}
\end{equation}
Because the objective depends only on $(\vect q,\tau,\vect\nu)$, the condensed and full-state forms share the same optimal control, time, buffer, and objective value for the current linearization. This equivalence does not extend to the nonlinear dynamics; nonlinear consistency is restored by the outer propagation.

\subsection{Structured Interval Solve}

With state order $(\vect r,\vect v,m)$, the state Jacobian has the approximate block form
\begin{equation}
F_x=\begin{bmatrix}
0&I&0\\
A_r&A_v&\vect a_m\\
0&0&0
\end{bmatrix}.
\end{equation}
For $c=\gamma\bar\tau$, solving $(I-cF_x)\vect w=\vect b$ yields
\begin{equation}
\vect w_r=\vect b_r+c\vect w_v,\qquad w_m=b_m,
\end{equation}
and leaves only
\begin{equation}
\underbrace{(I-cA_v-c^2A_r)}_{S_v}\vect w_v
=\vect b_v+cA_r\vect b_r+c\vect a_m b_m.
\label{eq:schur3}
\end{equation}
The implementation reuses a single pivoted $3\times3$ LU factorization of $S_v$ across 15 right-hand sides: seven state-transition columns, six endpoint-control columns, one time column, and one defect column.

The condensation requires every interval matrix $M_k=I-\gamma\bar\tau F_{x,k+1}$, equivalently its Schur complement $S_v$, to be nonsingular. The implementation rejects any factorization whose pivot falls below its numerical threshold; no rejection occurs in the nominal solve. An offline double-precision audit of the final 30 interval matrices yields a maximum two-norm condition number of \IntervalConditionMax\ and a minimum reciprocal one-norm condition estimate of \IntervalRcondMin. These values support the nonsingularity assumption for the reported trajectory but do not guarantee invertibility over a wider mission envelope.

At interval $k$, future controls cannot yet influence the current state, so the implementation propagates only the active prefix of $G_k$. For $N=31$, the number of propagated control-column instances is
\begin{equation}
\sum_{k=0}^{N-2}3(k+2)
=\frac{3}{2}(N-1)(N+2)=1485,
\end{equation}
instead of the 2790 instances required by an unconditional 93-column propagation---a 46.8\% reduction.

\subsection{Scaling and Rank}

We scale terminal position and velocity by \SI{500}{\meter} and
\SI{50}{\meter\per\second}, respectively, and then normalize each terminal row to unit Euclidean norm. The internal column scale is
\begin{equation}
D_c=\diag(\underbrace{1,\ldots,1}_{93},0.5,
\underbrace{0.25,\ldots,0.25}_{6}).
\end{equation}
With terminal physical scale $S_f$ and row normalization $R$, the solver receives
\begin{equation}
\widetilde H=RS_fHD_c,\qquad
\widetilde{\vect b}=RS_f\vect b.
\label{eq:scaledmatrix}
\end{equation}
The objective curvature and linear vector are transformed consistently. Because the six buffer columns are diagonal and nonzero for finite $\mu_\nu$, the condensed equality retains full row rank even when the early control sensitivities are nearly dependent.

\begin{table}[!t]
\centering
\caption{Dimensions of the convex subproblem for the default 31-node case}
\label{tab:dimensions}
\begin{tabular}{lrr}
\toprule
Quantity & Condensed xPIPG & Full-state ECOS \\
\midrule
Intermediate-state variables & 0 & 217 \\
Control variables & 93 & 93 \\
Final-time variables & 1 & 1 \\
Buffer and epigraph variables & 6 & 45 \\
Total primal variables & 100 & 356 \\
Equality constraints & 6 & 223 \\
Second-order cones & proximal & 64 \\
Dominant linear algebra & $6\times100$ products & sparse KKT factorization \\
\bottomrule
\end{tabular}
\end{table}

Complete condensation is most favorable when hard state-path constraints are absent or few. Condensing nodewise dynamic-pressure, heating, obstacle, or angle-of-attack constraints would introduce dense control rows, so such problems may instead favor partial condensation or the original sparse full-state form.

\section{Strongly Convex Surrogate and Latency-Oriented xPIPG}

\subsection{Separable Objective}

Let $\omega_1=\omega_N=1/2$ and all other trapezoidal weights equal one, and set $\Delta s=1/(N-1)$. Define
\begin{equation}
\varphi(\vect q)=\Delta s\sum_{k=1}^N\omega_k
\left[(1-\lambda)\norm{\vect q_k}_2
+\lambda\norm{\vect q_k}_2^2\right].
\end{equation}
The product $\tau\varphi(\vect q)$ is not jointly convex. The outer model therefore retains its complete convex dependence on the control and linearizes only the time factor:
\begin{equation}
\tau\varphi(\vect q)\approx
\bar\tau\varphi(\vect q)
+\varphi(\bar{\vect q})(\tau-\bar\tau).
\label{eq:objectiveapproximation}
\end{equation}
After condensation and scaling, the subproblem reads
\begin{equation}
\min_{\vect z}\ 
\frac{1}{2}\vect z^\mathsf TP\vect z+\vect c^\mathsf T\vect z
+g(\vect z)
\quad\mathrm{s.t.}\quad H\vect z=\vect b.
\label{eq:compositeproblem}
\end{equation}
The nonsmooth term is
\begin{equation}
g(\vect z)=\sum_{k=1}^Na_k\norm{\vect q_k}_2
+\sum_{k=1}^N\ind_{\mathcal U_k^p}(\vect q_k)
+\ind_{[\tau_{\min},\tau_{\max}]}(\tau),
\label{eq:nonsmoothterm}
\end{equation}
where
\begin{align}
P_{q,k}&=2\lambda\bar\tau\Delta s\omega_k+w_q,\\
\vect c_{q,k}&=-w_q\bar{\vect q}_k,\\
a_k&=(1-\lambda)\bar\tau\Delta s\omega_k,\\
P_\tau&=w_t,\qquad
c_\tau=\varphi(\bar{\vect q})-w_t\bar\tau.
\end{align}
Because the scaled buffer has unit curvature and no nonsmooth term, $P\succ0$, and the primal optimum of the mathematical problem in Eq.~\eqref{eq:compositeproblem} is unique. This uniqueness statement concerns the ideal surrogate; whether an iterative implementation reaches that optimizer depends on its actual node operator.

\subsection{Ideal-Surrogate Optimality Conditions and Implemented xPIPG Feedback}

If the exact joint proximal map of $g$ were evaluated, the optimality conditions for Eq.~\eqref{eq:compositeproblem} would be
\begin{subequations}\label{eq:kkt}
\begin{align}
\zero&\in P\vect z^\star+\vect c+\partial g(\vect z^\star)
+H^\mathsf T\vect\eta^\star,\\
\zero&=H\vect z^\star-\vect b.
\end{align}
\end{subequations}
These relations define the reference KKT point of the strongly convex surrogate. The timed implementation does not claim to solve them exactly: it replaces the exact control-block proximal operator by the one-pass map $\mathcal M_1$ defined below, while retaining exact clipping for time and the identity map for the unconstrained scaled buffer. The resulting latency-oriented xPIPG iteration is
\begin{subequations}\label{eq:pipg}
\begin{align}
\vect y^{j+1}
&=\mathcal M_1\!\left(
\vect\xi^j-\alpha(P\vect\xi^j+\vect c+H^\mathsf T\vect\eta^j)
\right),\\
\widehat{\vect\eta}^{j+1}
&=\vect\eta^j+\beta\left[H(2\vect y^{j+1}-\vect\xi^j)-\vect b\right],\\
\vect\xi^{j+1}
&=\vect\xi^j+\rho_e(\vect y^{j+1}-\vect\xi^j),\\
\vect\eta^{j+1}
&=\vect\eta^j+\rho_e(\widehat{\vect\eta}^{j+1}-\vect\eta^j),
\qquad \rho_e=1.55.
\end{align}
\end{subequations}
where $\rho_e=1.55$ makes the implemented method an extrapolated PIPG variant. Because $\mathcal M_1\ne\prox_{\alpha g}$ in general, Eq.~\eqref{eq:pipg} is an explicitly approximate fixed-point iteration. The KKT system in Eq.~\eqref{eq:kkt} is retained only to define the ideal convex reference, not as a theorem about the timed trajectory.

The step-size calculation requires an estimate of $\norm H_2^2=\lambda_{\max}(HH^\mathsf T)$. Because $H$ has six rows, the solver forms a $6\times6$ Gram matrix and applies 80 power iterations, giving $\widehat\sigma^2$. With $L=\max_jP_{jj}$, primal--dual ratio $r=1$, and empirical safety factor $s=0.97$,
\begin{equation}
\alpha=\frac{2s}
{L+\sqrt{L^2+4sr\widehat\sigma^2}},
\qquad \beta=r\alpha.
\label{eq:stepsize}
\end{equation}
Direct substitution gives
\begin{equation}
\alpha(L+\beta\widehat\sigma^2)=s<1.
\end{equation}
Finite power iteration returns a Rayleigh quotient, which is generally a lower estimate rather than a certified upper bound. Multiplication by $1.0000001$ only protects against the observed rounding scale; it does not prove $\widehat\sigma^2\ge\norm H_2^2$. For the audited final nominal subproblem, the relative estimate error is \PowerRelativeError, and the exact eigenvalue remains inside the 3\% empirical step margin. A Gershgorin bound computed offline is \GershgorinRatio\ times the exact eigenvalue. The present step sizes are therefore well resolved for the tested case, but they do not provide a universal stability certificate.

\subsection{Exact Component Operators and the One-Pass Node Map}

For the fuel term, the exact component operator is
\begin{equation}
\prox_{\alpha a_k\norm{\cdot}_2}(\vect u)=
\max\left(1-\frac{\alpha a_k}{\norm{\vect u}_2},0\right)\vect u.
\label{eq:groupshrink}
\end{equation}
This group shrinkage is well defined at zero and does not rotate thrust.

For the pointing cone, let
$a=\vect e^\mathsf T\vect u$,
$\vect u_\perp=\vect u-a\vect e$,
$b=\norm{\vect u_\perp}_2$, and
\begin{equation}
\vect w_\theta=\cos\theta\,\vect e
+\sin\theta\,\frac{\vect u_\perp}{b}.
\end{equation}
The cone projection is
\begin{equation}
\proj_{\mathcal K_\theta}(\vect u)=
\begin{cases}
\vect u, & b\le a\tan\theta,\\
\zero, & a\cos\theta+b\sin\theta\le0,\\
(a\cos\theta+b\sin\theta)\vect w_\theta, & \text{otherwise}.
\end{cases}
\label{eq:coneprojection}
\end{equation}
Projection onto the intersection of this cone and the unit ball is obtained by applying Eq.~\eqref{eq:coneprojection}, followed by radial clipping when necessary; positive radial scaling preserves cone membership. The lower-thrust half-space projection is
\begin{equation}
\proj_{\vect d^\mathsf T\vect q\ge q_{\min}}(\vect u)
=\vect u+\max(0,q_{\min}-\vect d^\mathsf T\vect u)\vect d.
\label{eq:halfspaceprojection}
\end{equation}

Let $\mathcal S_{\gamma}$ denote Eq.~\eqref{eq:groupshrink},
$\mathcal C=\mathcal K_\theta\cap\{\vect q:\norm{\vect q}_2\le1\}$,
and $\mathcal H_k=\{\vect q:(\vect d_k^p)^\mathsf T\vect q\ge q_{\min}\}$.
The timed node operator is exactly
\begin{equation}
\mathcal M_{1,k}(\vect u)=
\proj_{\mathcal H_k}\!\left[
\proj_{\mathcal C}\!\left(\mathcal S_{\alpha a_k}(\vect u)\right)
\right].
\label{eq:onepassmap}
\end{equation}
It performs one shrinkage, one cone--ball projection, and one half-space projection, with no persistent Dykstra corrections. In general,
\begin{equation}
\mathcal M_{1,k}\ne
\prox_{\alpha a_k\norm{\cdot}_2+\ind_{\mathcal C}+\ind_{\mathcal H_k}}.
\end{equation}
The distinction is structural. For example, choose unit vectors $\vect q$ and $\vect d$ on opposite boundaries of a $45^\circ$ cone such that $\vect q^\mathsf T\vect d=0$. Both lie in $\mathcal C$. A final half-space correction with $q_{\min}=0.4$ gives $\vect q'=\vect q+0.4\vect d$ and $\norm{\vect q'}_2=\sqrt{1.16}>1$, thereby breaking the unit-ball constraint.

An offline audit compares Eq.~\eqref{eq:onepassmap} against converged Dykstra splitting \cite{dykstra1983algorithm} over \ProxAuditSamples\ randomly generated three-dimensional inputs, with reference directions inside the pointing cone and thresholds spanning the nominal scale. The median and 99th-percentile map discrepancies are \ProxErrorMedian\ and \ProxErrorPnn, respectively; \ProxViolationFraction\ of the broad synthetic inputs have a set violation above $10^{-7}$. These results rule out interpreting one pass as an exact joint prox. The production choice is intentional: the map has fixed work and minimal memory traffic, and the complete trajectory is accepted only after independent nodewise physical checks. These checks safeguard the reported task-level output but do not restore exact KKT optimality of Eq.~\eqref{eq:compositeproblem}.

\subsection{Per-Update Cost and Stopping Criteria}

For $n_z=3N+7$ and terminal dimension $m_f=6$, the cost of one inner update is dominated by the products $H^\mathsf T\vect\eta$ and $H\vect y$, $N$ three-dimensional one-pass map evaluations, and the vector extrapolation:
\begin{equation}
\mathcal O(2m_fn_z+N)=\mathcal O(N).
\end{equation}
The spectral estimate is computed once per outer iteration. Persistent storage consists of the node trajectories, a $6\times100$ condensed matrix, a $93\times7$ terminal sensitivity, and several 100-vectors, all statically sized.

The first outer iteration uses inner tolerance $2\times10^{-3}$; subsequent tolerances are
\begin{equation}
\epsilon_{\mathrm{in}}=
\max\!\left\{2\times10^{-4},
\min(2\times10^{-3},0.05e_f)\right\}.
\end{equation}
Every fourth inner update checks
\begin{align}
\norm{H\vect y-\vect b}_\infty&\le\epsilon_{\mathrm{in}},\\
\frac{\norm{D_c(\vect y^j-\vect y^{j-1})}_2}
{\max(1,\norm{D_c\vect y^j}_2)}
&\le\epsilon_{\mathrm{in}}.
\end{align}
This check interval avoids copying a 100-vector and evaluating two norms on the remaining three iterations. The minimum and maximum inner budgets are 20 and 6000 updates.

\section{Fixed-Array Implementation}

\subsection{Data Layout and Hot Loop}

The C17 solver fixes $N=31$, a state dimension of seven, a control dimension of 93, a primal dimension of 100, and a dual dimension of six at compile time. The primal ordering is
\begin{equation}
[q_{1x},q_{1y},q_{1z},\ldots,q_{Nz},
\tau,\nu_1,\ldots,\nu_6].
\end{equation}
The three control components of each node are stored contiguously, which enables a fused nodewise map. The row-major $H[6][100]$ matrix occupies approximately \SI{2.4}{\kilo\byte} in single precision, and every component of $H^\mathsf T\vect\eta$ is a fixed six-term dot product. No memory is dynamically allocated inside the outer approximation or the xPIPG loop.

Before the hot loop, the solver caches
\begin{equation}
p_j=1-\alpha P_{jj},\qquad
d_j=-\alpha c_j,\qquad
\widehat H=\alpha H,
\end{equation}
so the primal trial component is
\begin{equation}
u_j=p_j\xi_j+d_j-\sum_{i=1}^6\widehat H_{ij}\eta_i.
\end{equation}
Only one explicit $\widehat H\vect y$ product is evaluated per iteration; the second follows from
\begin{equation}
\widehat H\vect\xi^{j+1}
=(1-\rho_e)\widehat H\vect\xi^j
+\rho_e\widehat H\vect y^{j+1}.
\end{equation}

\subsection{End-to-End Solve Workflow}

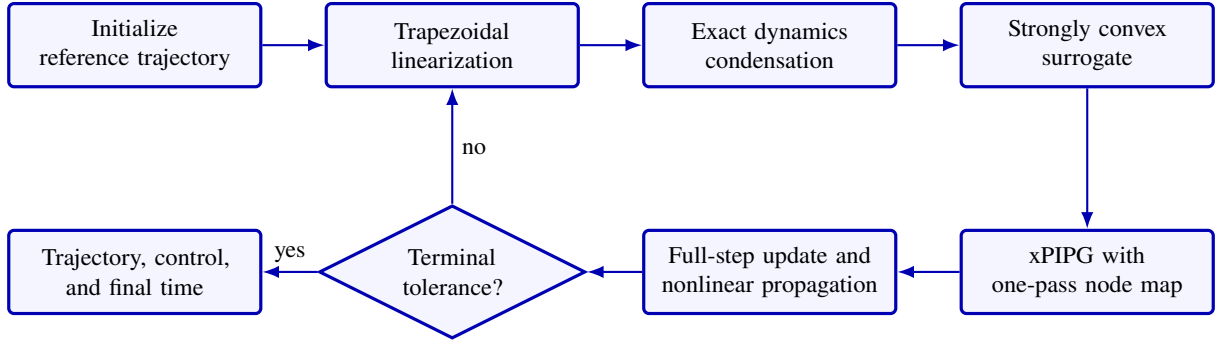
\begin{figure*}[!t]
\centering
\begin{tikzpicture}[
font=\small,
main/.style={draw=blue!70!black,very thick,rounded corners=2pt,
align=center,minimum height=11mm,text width=31mm,fill=blue!4},
decision/.style={draw=blue!70!black,very thick,diamond,aspect=2.0,
align=center,inner sep=1.5pt,text width=20mm,fill=blue!4},
arrow/.style={-{Latex[length=2.2mm]},thick,draw=blue!70!black}
]
\node[main] at (0,1.8) (init) {Initialize\\reference trajectory};
\node[main] at (4.2,1.8) (lin) {Trapezoidal\\linearization};
\node[main] at (8.4,1.8) (con) {Exact dynamics\\condensation};
\node[main] at (12.6,1.8) (sub) {Strongly convex\\surrogate};
\node[main] at (12.6,-1.2) (pipg) {xPIPG with\\one-pass node map};
\node[main] at (8.4,-1.2) (prop) {Full-step update and\\nonlinear propagation};
\node[decision] at (4.2,-1.2) (stop) {Terminal\\tolerance?};
\node[main] at (0,-1.2) (out) {Trajectory, control,\\and final time};

\draw[arrow] (init)--(lin);
\draw[arrow] (lin)--(con);
\draw[arrow] (con)--(sub);
\draw[arrow] (sub)--(pipg);
\draw[arrow] (pipg)--(prop);
\draw[arrow] (prop)--(stop);
\draw[arrow] (stop)--node[above]{yes}(out);
\draw[arrow] (stop.north) -- node[right,pos=0.48]{no} (lin.south);
\end{tikzpicture}
\caption{Online workflow of the condensed sequential convex approximation. Offline reference solvers and validation routines lie outside the timed path.}
\label{fig:workflow}
\end{figure*}

Figure~\ref{fig:workflow} shows only the online sequence. Each outer update performs one linearization, one exact condensation, one strongly convex surrogate construction, one latency-oriented xPIPG pass up to the stopping criterion, and one nonlinear propagation. The ``no'' branch returns directly from the decision node to the next linearization. MATLAB scripting, high-order integration, the converged Dykstra reference, and full-state ECOS serve as offline audits and are deliberately omitted from the online flow.

The timed C build uses Visual Studio 2022 Release x64 with C17,
\texttt{/O2}, \texttt{/Oi}, \texttt{/Ot}, \texttt{/GL},
\texttt{/arch:AVX2}, and \texttt{/fp:fast}. These settings expose the implementation potential on the test processor; they do not establish flight qualification. Strict floating-point agreement, worst-case execution time on the target processor, cache and frequency controls, timeout behavior, and retention of a previously validated trajectory all remain necessary for flight software.

\begin{table*}[!t]
\centering
\caption{Observed per-stage workload for one nominal condensed C solve}
\label{tab:workload}
\begin{tabular}{@{}
>{\raggedright\arraybackslash}p{0.19\textwidth}
>{\raggedright\arraybackslash}p{0.12\textwidth}
>{\raggedright\arraybackslash}p{0.26\textwidth}
>{\raggedright\arraybackslash}p{0.30\textwidth}@{}}
\toprule
Stage & Count & Principal operation & Fixed structure or bound \\
\midrule
Node linearization & $4\times31=124$ nodes & baseline dynamics, three velocity differences, and the analytic control Jacobian & each unique node is evaluated once per outer iteration \\
Interval model and condensation & $4\times30=120$ intervals & one $3\times3$ LU factorization, 15 right-hand sides, active-column propagation & a fixed set of 30 intervals and no candidate steps \\
Spectral estimate and step size & 4 & 21 length-100 dot products and 80 six-dimensional power iterations & once per surrogate; an empirical rather than a certified upper bound \\
xPIPG hot loop & 336 & terminal matrix products and 31 three-dimensional one-pass maps & static arrays; the stopping test runs every fourth update \\
Nonlinear propagation & $4\times30=120$ intervals & implicit trapezoidal Newton solve & at most 12 Newton steps per interval \\
Final physical verification & 1 & thrust, pointing, mass, impulse, time, and terminal checks & independent of dual variables and scaling \\
\bottomrule
\end{tabular}
\end{table*}

\section{Numerical Results}

\subsection{Case Definition and Timing Protocol}

Table~\ref{tab:caseparameters} lists the nominal case. The initial velocity has a magnitude of \SI{210}{\meter\per\second} and points approximately toward the target. The final time is initialized at the midpoint of the allowable interval. The first state reference is a bounded endpoint interpolation; after each convex solve, the original implicit nonlinear dynamics generate the next reference.

\begin{table}[!t]
\centering
\caption{Nominal powered-landing parameters}
\label{tab:caseparameters}
\begin{tabular}{lll}
\toprule
Group & Parameter & Value \\
\midrule
Initial state & $\vect r_0$ (\si{\meter}) & $[500,2500,500]$ \\
& $\vect v_0$ (\si{\meter\per\second}) & $[-49.5026,-197.9873,-49.5026]$ \\
Mass & $m_0,m_{\mathrm{dry}}$ (\si{\kilogram}) & $32000,\ 25600$ \\
Engine & $I_{sp}$ (\si{\second}) & $282$ \\
& $T_{\min},T_{\max}$ (\si{\kilo\newton}) & $338,\ 845$ \\
Geometry & $S$ (\si{\meter\squared}) & $10.752$ \\
& $\theta_{\max}$ & $45^\circ$ \\
Time and grid & $[t_{f,\min},t_{f,\max}]$ (\si{\second}) & $[15,45]$ \\
& $N$ & 31 \\
Atmosphere & $\rho_0$ (\si{\kilogram\per\meter\cubed}) & 1.225 \\
& $H_s$ (\si{\meter}), $C_D$ & $8500,\ 1.5$ \\
Objective & $\lambda,w_q,w_t$ & $0.01,\ 0.02,\ 1.0$ \\
Tolerance & $e_f,\epsilon_{\mathrm{in,min}}$ & $10^{-3},\ 2\times10^{-4}$ \\
\bottomrule
\end{tabular}
\end{table}

The test platform is an Intel Core i7-10875H running Windows 11 and Visual Studio 2022. A high-resolution performance counter brackets the complete sequence of initialization, outer approximation, inner solve, nonlinear propagation, and nodewise physical verification; process startup, console output, file writing, and plotting are excluded. The nominal tests use 20 warmups followed by 2,000 measured solves for xPIPG, and five warmups followed by 200 measured solves for ECOS. The Monte Carlo cases use five repetitions per sample for xPIPG and three for ECOS, retaining the within-sample median. The grid tests use 300 repetitions for xPIPG and 20 for ECOS. Both executables share the same task-level success condition: a six-dimensional nonlinear terminal residual no larger than $10^{-3}$, together with final-time, nodewise thrust-magnitude, pointing-cone, dry-mass, and mass--impulse checks. Because the ECOS P99 is inferred from only 200 timings, it is a descriptive tail statistic rather than a stable extreme-quantile estimate.

\subsection{Nominal Solution and Microsecond-Class Latency}

\begin{figure*}[!t]
\centering
\includegraphics[width=0.98\textwidth]{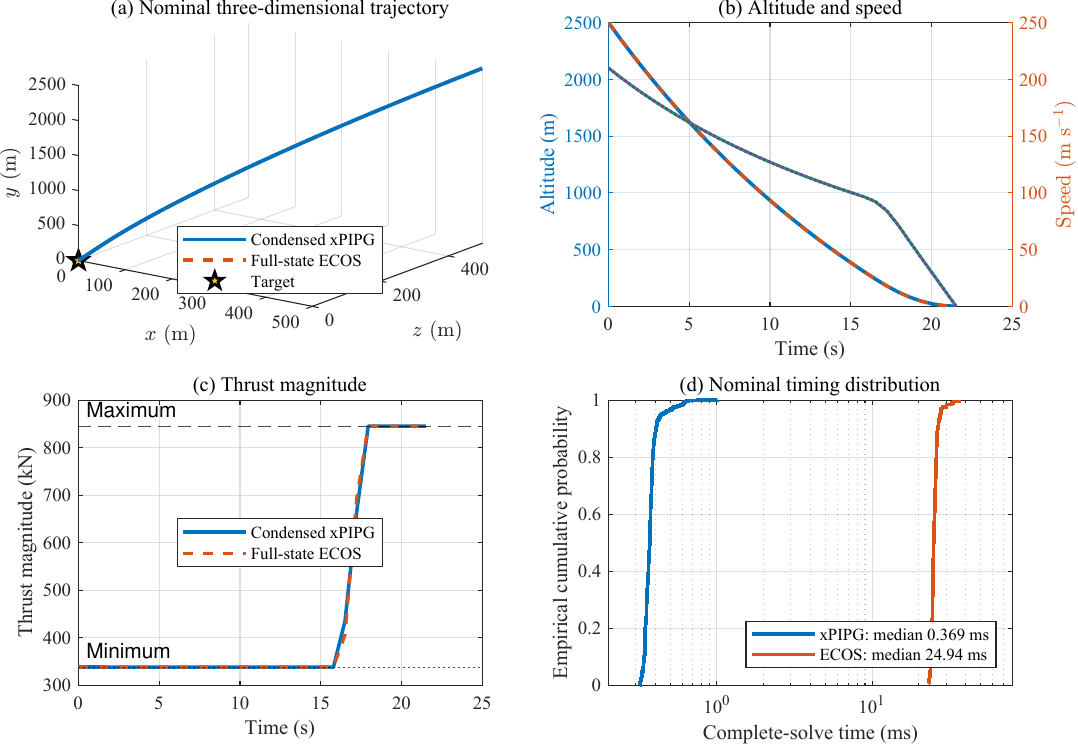}
\caption{Nominal trajectory, states, thrust history, and end-to-end solve-time distributions. All labels and annotations are regenerated in English from the archived benchmark data.}
\label{fig:nominal}
\end{figure*}

Figure~\ref{fig:nominal} compares the nominal trajectories and the end-to-end solve-time distributions. Both solutions reach the target at zero terminal velocity and satisfy the lower and upper thrust-magnitude bounds of \SIrange{338}{845}{\kilo\newton} and the $45^\circ$ pointing cone.

\begin{table}[!t]
\centering
\caption{Nominal solution and end-to-end solve-time statistics}
\label{tab:nominal}
\begin{tabular}{lrrrrrr}
\toprule
Method & Outer & Inner & Median & P99 & $t_f$ & Propellant \\
& & & (ms) & (ms) & (s) & (kg) \\
\midrule
Condensed xPIPG & 4 & 336 & 0.374 & 0.512 & 21.5283 & 3464.57 \\
Full-state ECOS & 4 & 64 & 24.96 & 32.55 & 21.5252 & 3463.61 \\
\bottomrule
\end{tabular}
\end{table}

The condensed solver terminates with $e_f=8.77\times10^{-4}$; ECOS reaches $9.91\times10^{-6}$. The final times differ by only \SI{2.79}{\milli\second}, and the propellant consumption by \SI{0.96}{\kilogram}. The maximum nodewise position and velocity differences are \SI{0.435}{\meter} and \SI{0.371}{\meter\per\second}, respectively. These small differences support physical agreement between the two solutions, but they are not an optimality-error claim, because the objective and buffer norms differ.

The xPIPG median of \SI{0.374}{\milli\second} corresponds to \SI{374}{\micro\second}, with a P99 of \SI{512}{\micro\second}. The median end-to-end speed ratio relative to the present ECOS implementation is 66.8. Although xPIPG performs more inner updates, each update involves only fixed-size matrix--vector products and nodewise one-pass node maps, whereas an ECOS interior-point step updates the cone scaling and factorizes a full-state KKT system. The ratio therefore compares two complete formulations and implementations; it does not isolate the benefit of xPIPG or of condensation alone. More importantly, the absolute \SI{374}{\micro\second} result places the implementation in the submillisecond regime. Within the scope surveyed in Table~\ref{tab:literatureposition}, it is the fastest sequential-convex powered-landing solve known to the authors and the only surveyed result below \SI{1}{\milli\second}.

\subsection{Independent Continuous-Time Audit}

The online propagation uses the same implicit trapezoidal transcription as the outer model, so agreement at the nodes alone cannot exclude intersample error. A separate offline audit therefore propagates the archived C control with first-order hold using MATLAB \texttt{ode113} at a relative tolerance of $10^{-11}$ and 101 query points per interval. This audit changes both the integration method and the intersample evaluation density. Table~\ref{tab:continuousaudit} compares the terminal state and the dense constraint checks.

\begin{table}[!t]
\centering
\caption{Independent first-order-hold continuous-time audit}
\label{tab:continuousaudit}
\begin{tabular}{lrr}
\toprule
Quantity & Node/trapezoidal & FOH/\texttt{ode113} \\
\midrule
Terminal position error (\si{\meter}) & 0.439 & \ContinuousPositionError \\
Terminal velocity error (\si{\meter\per\second}) & 0.040 & \ContinuousVelocityError \\
Normalized terminal infinity norm & $8.77\times10^{-4}$ & \ContinuousTerminalResidual \\
Minimum-thrust violation (\si{\newton}) & 0 & \ContinuousMinimumThrustViolation \\
Maximum-thrust violation (\si{\newton}) & 0 & 0.048 \\
Pointing violation (\si{\newton}) & 0 & \ContinuousMaximumPointingViolation \\
\bottomrule
\end{tabular}
\end{table}

The high-order propagation reaches the target with smaller terminal errors than the node trajectory, but dense sampling detects a \SI{1.688}{\newton} lower-thrust undershoot, approximately $5.0\times10^{-6}$ of $T_{\min}$. This undershoot is negligible for the present numerical demonstration, yet it is a genuine intersample violation under a zero-tolerance interpretation. Flight use would require tightening the node lower bound, enforcing a continuous-time margin, or changing the hold parameterization. The audit therefore supports the continuous-dynamics accuracy of the nominal result while preventing nodewise validation from being presented as an exact guarantee of continuous-time feasibility.

\subsection{Common-Sample Monte Carlo Results}

Each component of the initial position and velocity is perturbed independently by a uniform relative factor between $-10$ and $10\%$. A fixed random seed generates 100 six-dimensional samples, which are written once and read by both executables. Every case starts from the default initialization, with no sample-specific parameter tuning.

\begin{figure*}[!t]
\centering
\includegraphics[width=0.98\textwidth]{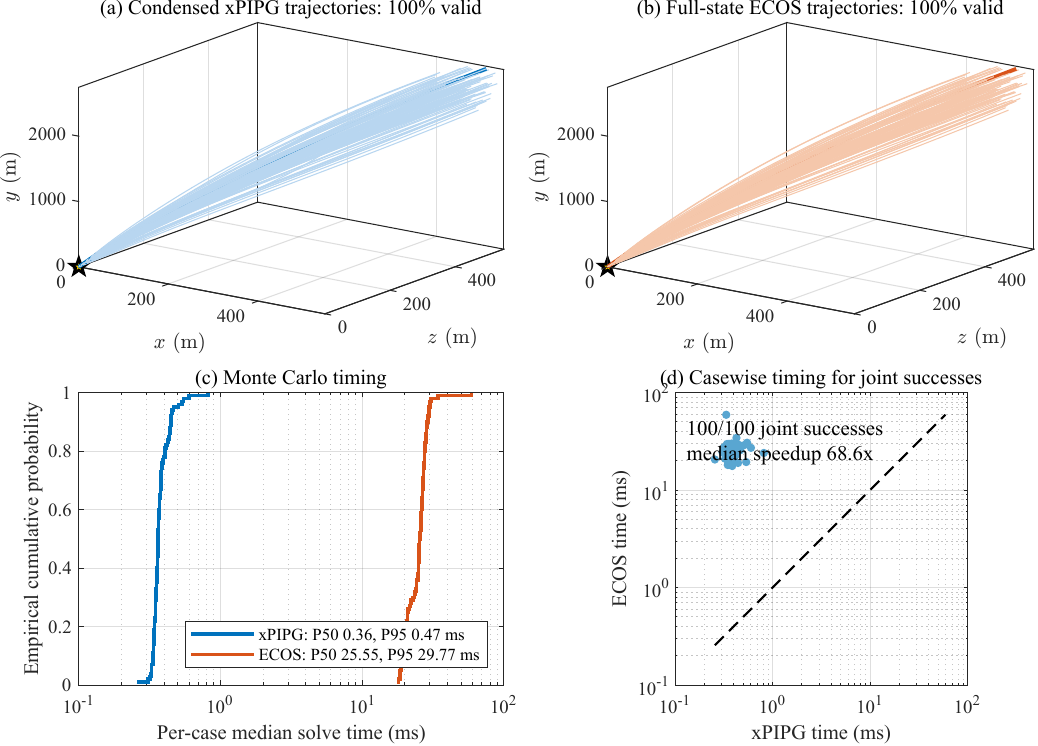}
\caption{Trajectory families and timing distributions for 100 common initial-state perturbations. Failed cases would appear as dashed red trajectories; none occur in this test.}
\label{fig:montecarlo}
\end{figure*}

\begin{table}[!t]
\centering
\caption{Statistics over 100 common perturbed initial states}
\label{tab:montecarlo}
\begin{tabular}{lrrrrrr}
\toprule
Method & Valid & P50 & P95 & P99 & Maximum & Mean outer \\
& (\%) & (ms) & (ms) & (ms) & (ms) & \\
\midrule
Condensed xPIPG & 100 & 0.36 & 0.47 & 0.60 & 0.81 & 4.25 \\
Full-state ECOS & 100 & 25.55 & 29.77 & 34.56 & 59.19 & 3.73 \\
\bottomrule
\end{tabular}
\end{table}

Both implementations validate all 100 cases. Among the cases solved by both methods, the casewise median ECOS/xPIPG speed ratio is 68.6. xPIPG uses slightly more outer iterations on average, but the small condensed inner model keeps its P95 below \SI{0.5}{\milli\second}. For both solvers, the maxima exceed P99, which illustrates why a statistical desktop tail must not be presented as a target-processor worst-case execution time. The $\pm10\%$ sample is a local engineering screen, not evidence of robustness over a mission-wide envelope or against aerodynamic, wind, mass, and actuator uncertainty.

\subsection{Grid-Size Scaling}

The scan uses five grids from 15 to 51 nodes while holding the physical parameters, stopping criteria, and algorithm parameters fixed. All ten combinations of solver and grid pass independent physical validation.

\begin{figure*}[!t]
\centering
\includegraphics[width=0.96\textwidth]{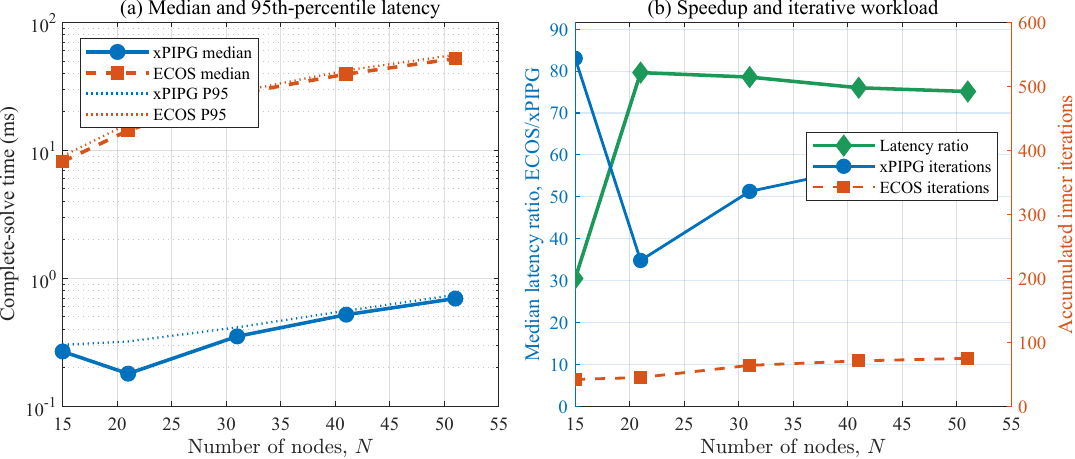}
\caption{End-to-end solve time and accumulated iterative workload as the number of discretization nodes varies.}
\label{fig:gridscaling}
\end{figure*}

\begin{table*}[!t]
\centering
\caption{End-to-end scaling with the number of discretization nodes}
\label{tab:gridscaling}
\begin{tabular}{rrrrrrrr}
\toprule
$N$ & xPIPG median & xPIPG P95 & ECOS median & ECOS P95 &
Median ratio & xPIPG inner & ECOS inner \\
& (ms) & (ms) & (ms) & (ms) & & & \\
\midrule
15 & 0.268 & 0.303 & 8.18 & 9.04 & 30.5 & 544 & 42 \\
21 & 0.180 & 0.321 & 14.33 & 16.14 & 79.6 & 228 & 45 \\
31 & 0.352 & 0.414 & 27.66 & 28.66 & 78.6 & 336 & 64 \\
41 & 0.520 & 0.558 & 39.55 & 42.23 & 76.0 & 372 & 71 \\
51 & 0.696 & 0.744 & 52.27 & 55.93 & 75.1 & 408 & 75 \\
\bottomrule
\end{tabular}
\end{table*}

The condensed median latency ranges from \SI{0.180}{\milli\second} to
\SI{0.696}{\milli\second}; the ECOS median ranges from \SI{8.18}{\milli\second} to
\SI{52.27}{\milli\second}. The runtime is not strictly monotone across the two coarsest grids, because refining the grid also changes the trapezoidal-model error, the conditioning of the terminal sensitivity, and the location of thrust-set transitions. The 15-node case therefore requires 544 xPIPG updates, compared with 228 for 21 nodes. From 21 to 51 nodes, the median speed ratio remains between 75.1 and 79.6.

\begin{table}[!t]
\centering
\caption{Physical outputs of the xPIPG solver under grid refinement}
\label{tab:gridphysics}
\begin{tabular}{rrrr}
\toprule
$N$ & $t_f$ (\si{\second}) & Propellant (\si{\kilogram}) & $e_f$ \\
\midrule
15 & 21.5884 & 3469.33 & $2.25\times10^{-4}$ \\
21 & 21.5456 & 3465.27 & $5.54\times10^{-4}$ \\
31 & 21.5283 & 3464.57 & $8.77\times10^{-4}$ \\
41 & 21.5341 & 3465.72 & $9.89\times10^{-4}$ \\
51 & 21.5445 & 3466.69 & $4.41\times10^{-4}$ \\
\bottomrule
\end{tabular}
\end{table}

Across 15--51 nodes, the final time varies by \SI{0.061}{\second} and the propellant consumption by \SI{4.76}{\kilogram}, or 0.14\% of the nominal value. The nonmonotone propellant sequence shows that this scan indicates discretization consistency rather than proving convergence to a continuous optimum.

\section{Discussion}

\subsection{Why the Runtime Reaches Hundreds of Microseconds}

The measured latency follows from structure at four interacting levels. At the model level, exact dynamics condensation eliminates the entire intermediate-state chain before optimization begins. At the algorithm level, only six terminal dual variables remain, and every control operation is three-dimensional. At the numerical level, row and column scaling together with a small Gram-matrix estimate remove any need for KKT factorization. At the implementation level, static single-precision arrays, shared node Jacobians, $3\times3$ Schur complements, active sensitivity columns, a one-pass node map, product recurrences, and downsampled stopping checks cut memory traffic and branching. No single transformation accounts for the \SI{374}{\micro\second} median.

The timing comparison is deliberately end to end: the xPIPG and ECOS bars combine different state representations, objective details, buffer norms, inner tolerances, arithmetic modes, memory layouts, warmup policies, and repetition counts. The results therefore support the computational value of the complete condensed design; they do not isolate the speedup attributable to xPIPG, to condensation, or to compiler choices. In particular, the factor 66.8 is not the outcome of a same-cone solver competition, and the literature comparison in Table~\ref{tab:literatureposition} is likewise contextual. The defensible efficiency claim is the absolute measured latency under the stated timing boundary, not a universal solver ranking.

\subsection{Strong Convexity and the Implemented Approximation}

The $\lambda=0.01$ energy fraction keeps the Euclidean fuel norm dominant while contributing smooth curvature. The fixed $w_q$ and $w_t$ weights further regularize the change between adjacent outer iterates. Because these choices modify every local subproblem, the computed solution is not a mathematically exact pure-fuel optimum. Reducing the curvature moves the solution toward bang--bang fuel behavior but can slow or destabilize a first-order fixed point near shifting active sets; increasing it smooths the thrust profile but may cost additional impulse.

Strong convexity guarantees a unique primal solution only for the ideal condensed surrogate of Eq.~\eqref{eq:compositeproblem}. Because the timed one-pass node map is not its exact joint prox, that uniqueness does not imply that the implemented iteration returns the unique optimizer. Likewise, the 80-step power estimate is not a certified spectral upper bound. The random proximal audit and the exact $6\times6$ eigenvalue audit quantify these two approximations rather than concealing them. The evidence supports a fast, physically screened approximate solution; it does not establish full KKT convergence.

\subsection{Outer Update and Failure Handling}

The outer method is best described as a proximal-regularized sequential convex approximation. It accepts every computed update in full, with no actual-to-predicted reduction ratio, no backtracking, no candidate rejection, and no dynamic trust-region radius. Although the fixed $w_q$ and $w_t$ terms discourage large changes, they are not equivalent to a trust-region acceptance mechanism. Conventional convergence statements for safeguarded SCP therefore do not transfer directly.

The present initialization and all 100 local perturbations converge, but failures are expected for sufficiently remote initial guesses, ill-conditioned interval matrices, inconsistent final-time guesses, abrupt changes in the minimum-thrust active set, or one-pass maps that fail the final physical screen. The implementation responds deterministically: it rejects the candidate trajectory and retains a previously validated command. The current experiments map neither these failure boundaries nor the sensitivity to all four regularization and penalty parameters.

\subsection{Consolidated Validity Boundary}

The favorable condensed dimension rests on a small terminal constraint set and few hard state-path constraints. Dense nodewise limits on dynamic pressure, heating, obstacles, glide slope, or angle of attack can make complete condensation unattractive; partial condensation or a full-state sparse formulation should then be considered.

The validated claim covers only a three-degree-of-freedom point-mass model with constant $C_D$, no hard state-path constraints, a 31-node nominal grid, a local $\pm10\%$ initial-state sample, and a desktop i7 processor. The model omits attitude and angular-rate dynamics, gimbal-rate limits, structural loads, terrain, wind and parameter uncertainty, and plume--surface interaction. Node checks do not imply exact intersample feasibility, as the small FOH lower-thrust undershoot demonstrates. The statistical desktop tail is not a flight-processor worst-case execution time.

Flight application further requires a vehicle-specific aerodynamic database, six-degree-of-freedom coupling where appropriate, explicit continuous-time margins, mission-envelope Monte Carlo campaigns that include failures and infeasible cases, identical-precision ablations, processor- and hardware-in-the-loop testing, a certified timeout policy, and worst-case timing on the intended processor. These items are requirements for future extension, not claims of the present paper.

\input{expanded_validation_and_reproducibility.tex}

\section{Conclusions}

This paper has presented a latency-oriented sequential convex approximation for powered-landing optimization. For each local linearized model, exact dynamics condensation replaces 210 interval equations with six terminal equalities. A fuel--energy objective with fixed quadratic proximal terms gives the ideal surrogate its strong-convexity curvature. The xPIPG inner solver and fixed-size arrays then reduce every online update to small matrix--vector products and three-dimensional node operations.

The nominal C implementation completes one plan in four outer and 336 inner updates. Its median end-to-end solve time is \SI{374}{\micro\second}, with a P99 of \SI{512}{\micro\second}. All 100 common $\pm10\%$ initial-state perturbations pass task-level validation, and the median stays below \SI{0.7}{\milli\second} through 51 nodes. Independent first-order-hold integration yields a terminal position error of \SI{0.183}{\meter} and reveals a small \SI{1.688}{\newton} intersample lower-thrust violation.

The result should therefore be read as a low-latency, physically screened approximate solution rather than as an exact KKT solution of the ideal strongly convex surrogate. Both the one-pass node map and the finite power iteration are deliberate low-cost approximations. Future work will pursue exact low-cost node active sets, certified spectral bounds, safeguarded outer updates, partial condensation, wider mission envelopes, and representative flight processors.

\section*{Funding Sources}

This research received no external funding.

\section*{Declarations}

None.

\input{appendix_algorithms.tex}

\bibliography{references}

\end{document}

%% file: review_audit/reviewer_audit_values.tex
\newcommand{\ProxAuditSamples}{50000}
\newcommand{\ProxErrorMedian}{3.022e-03}
\newcommand{\ProxErrorPnn}{5.101e-01}

\newcommand{\ProxViolationFraction}{9.51\%}

\newcommand{\PowerRelativeError}{1.000e-07}
\newcommand{\GershgorinRatio}{1.1325}
\newcommand{\IntervalConditionMax}{1.436}
\newcommand{\IntervalRcondMin}{5.333e-01}
\newcommand{\ContinuousPositionError}{0.183}
\newcommand{\ContinuousVelocityError}{0.011}
\newcommand{\ContinuousTerminalResidual}{3.606e-04}
\newcommand{\ContinuousMinimumThrustViolation}{1.688}
\newcommand{\ContinuousMaximumPointingViolation}{0.000}

%% file: expanded_validation_and_reproducibility.tex
\section{Parameter Roles, Validation Layers, and Reproducibility}

\subsection{Layered Effects of the Algorithm Parameters}

These parameters do not act at a common level, and tuning them as if they did can be misleading. The fuel--energy weight $\lambda$ changes both the physical preference and the control curvature of each local objective. The fixed proximal weights $w_q$ and $w_t$ limit how far the control and the final time may move between consecutive linearizations. The virtual-buffer penalty $\mu_\nu$ governs how readily the linearized terminal equality may borrow artificial feasibility. The quantities $r$, $s$, $\rho_e$, and the stopping tolerance primarily control the fixed-point path and the work required for a given surrogate. Table~\ref{tab:parameterroles} separates these effects. The distinction matters because varying $\lambda$, $w_q$, or $w_t$ changes the surrogate itself, whereas varying an admissible xPIPG step size or extrapolation factor leaves the optimizer of the ideal problem with an exact proximal map unchanged.

\begin{table*}[t]
\caption{Roles and mismatch signatures of the principal parameters}
\label{tab:parameterroles}
\centering
\small
\setlength{\tabcolsep}{4pt}
\begin{tabularx}{\textwidth}{p{0.08\textwidth}p{0.19\textwidth}X X p{0.12\textwidth}}
\toprule
Parameter & Direct role & Typical symptom when too small & Typical symptom when too large & Present value\\
\midrule
$\lambda$ & quadratic thrust curvature & near-switching control with sensitive active-set changes & excessive smoothing and possibly larger impulse & 0.01\\
$w_q$ & proximal regularization of the control & large control changes across outer iterations, with model mismatch & reference-control locking and slow outer progress & 0.02\\
$w_t$ & proximal regularization of the final time & rapid time changes amplify the bilinear error & final time remains close to its initial value & 1.0\\
$\mu_\nu$ & terminal virtual buffer & easy surrogate feasibility but large nonlinear terminal error & stiff scaling while the early reference is still inaccurate & starts at 50 and increases\\
$s$ & primal--dual stability margin & conservative steps and more inner work & insufficient margin under spectral underestimation & 0.97\\
$\rho_e$ & synchronous extrapolation & reverts toward the slower basic-PIPG path & oscillation near active-set changes & 1.55\\
$\epsilon_{\rm in}$ & inner stopping accuracy & wasted work on an inaccurate early model when overly strict & an error floor limits the outer terminal accuracy & $2{\times}10^{-3}$ to $2{\times}10^{-4}$\\
\bottomrule
\end{tabularx}
\end{table*}

The outer loop accepts each computed candidate in full and never evaluates an actual-to-predicted reduction ratio. Accordingly, $w_q$ and $w_t$ act as a fixed soft proximal regularization rather than as a trust region with radius updates and candidate rejection. This choice yields a single execution path per outer update, with no repeated candidate construction or propagation; the price is the absence of automatic recovery from severe linearization mismatch. The 100 local perturbations show that the selected parameters are adequate in the tested neighborhood; they do not establish convergence from arbitrary initial conditions or under arbitrary aerodynamic models.

\subsection{Three Distinct Levels of Correctness}

The method is best interpreted through three separate statements, which should not be conflated.

\begin{enumerate}
\item \textbf{Condensation level.} If every $M_k$ and every scaling transformation is nonsingular, the state recursion amounts to block Gaussian elimination of the full-state linear equality system. The condensed and uncondensed mathematical linearizations then share the same control--time feasible projection, and every condensed point has a unique recovered state sequence. The condensation at this level is exact.
\item \textbf{Surrogate level.} The small energy term and the fixed quadratic proximal terms together give the smooth part of the ideal condensed surrogate positive curvature. With an exact joint proximal operator and admissible steps, that ideal problem has a unique primal solution characterized by its KKT system.
\item \textbf{Implementation level.} The timed code replaces the joint proximal operator with a one-pass node map and obtains its step from finite power iteration. Its hot loop is therefore a latency-oriented approximate fixed-point iteration. Exact KKT convergence of the ideal surrogate does not transfer to this implementation; its output is instead screened through nonlinear terminal residuals, nodewise physics, and an independent continuous-time propagation.
\end{enumerate}

This separation also explains why ``the final trajectory passes the physical checks'' and ``the inner loop returns the exact strongly convex surrogate optimizer'' are different propositions: the first supports the engineering latency result reported here, whereas the second would require an exact joint proximal map, a certified spectral bound, and complete KKT residuals---none of which is claimed for the timed configuration.

\subsection{End-to-End Work and Resident Storage}

For $N=31$, one nominal plan contains 124 unique node linearizations, 120 interval models with their $3\times3$ Schur factorizations, four $6\times6$ Gram-matrix estimates, 336 inner updates, and four nonlinear propagations. At interval $k$, active-column propagation touches only the first $3(k+1)$ control columns; future columns remain structurally zero. Compared with propagating all 93 columns at every interval, this removes approximately half of the early dense products and avoids forming a global $210\times217$ state-chain matrix.

The principal resident arrays are the $31\times7$ state trajectory, the $31\times3$ control trajectory, the $6\times100$ condensed matrix, the terminal sensitivities, and several vectors of length 100. In row-major single precision, $H$ occupies approximately \SI{2.4}{\kilo\byte}, and every entry of $H^\mathsf T\vect\eta$ is a fixed six-term dot product. The three control components are stored contiguously, so the shrinkage, cone--ball projection, and minimum-thrust half-space projection fuse into a single node loop. The stopping test runs only every fourth update, so the other three updates neither copy the previous 100-vector nor evaluate two norms. Together, these details explain how 336 inner updates remain compatible with an end-to-end solve time in the hundreds of microseconds.

\subsection{Reproducibility Protocol and Result-Freezing Rule}

The top-level reproduction entry point is \path{run_paper_tests.ps1} under \path{test}. It calls \path{build_and_run_core_experiments.ps1} to generate the repeated nominal timings, the 100 common Monte Carlo cases, and the five grid-size cases. The script \path{generate_core_experiment_results.m} then rebuilds the statistics, tables, and English figures from the archived CSV files. The proximal, spectral, interval-conditioning, and continuous-time checks are offline audits and are therefore excluded from the reported online timing.

A reproduction run must freeze the physical parameters, initial condition, random seed, grid, compiler options, floating-point mode, warm-up count, repeat count, stopping conditions, and timing boundary. Changing any one of these items requires regenerating all tables and figures; individual numbers must not be substituted in isolation. Verification proceeds at four levels: agreement between the condensed prediction and the full linearized dynamics chain; inner equality and fixed-point residuals; terminal error under the original nonlinear trapezoidal propagation; and independent FOH high-order integration with intersample control checks. A case counts as successful only when all task-level online checks pass.

%% file: appendix_algorithms.tex
\appendix

\section{Supplementary Analyses of the Node Map and Spectral Estimate}

\subsection{Geometry and Counterexample for the One-Pass Node Map}

The ideal nonsmooth node term combines the fuel norm, the pointing cone, the unit ball, and the local minimum-thrust half-space. Group shrinkage,
\begin{equation}
\mathcal S_\kappa(\vect z)=
\max\left(1-\frac{\kappa}{\norm{\vect z}_2},0\right)\vect z ,
\end{equation}
is the exact proximal operator of the fuel norm. The cone--ball and half-space components also have low-dimensional closed-form projections. Applying these maps in series, however, generally yields neither the joint proximal operator of their sum nor the projection onto the full intersection.

A direct counterexample shows how the last half-space step can destroy unit-ball feasibility. Choose unit vectors $\vect q$ and $\vect d$ on opposite boundaries of a 45-deg pointing cone such that $\vect q^\mathsf T\vect d=0$. Both vectors satisfy the cone and unit-ball constraints. For $q_{\min}=0.4$, the final minimum-thrust half-space projection can produce
\begin{equation}
\vect q'=\vect q+0.4\vect d,\qquad
\norm{\vect q'}_2=\sqrt{1.16}>1.
\end{equation}
The production node operation must therefore be understood as a fixed-work approximation, and the 50,000-point offline Dykstra audit quantifies the resulting discrepancy. Final physical screening can reject an infeasible trajectory, but it cannot retroactively prove that the inner loop satisfies the ideal KKT system.

\subsection{Spectral Estimate, Certified Bounds, and the Latency Choice}

Finite power iteration on $G=HH^\mathsf T\succeq0$ returns a Rayleigh quotient that generally approaches $\lambda_{\max}(G)$ from below. Multiplication by $1.0000001$ therefore does not yield a rigorous upper bound. In the final nominal subproblem, the 80-step estimate is 1.580279627, whereas a double-precision symmetric eigensolve gives 1.580279469. The small positive error demonstrates adequate resolution only for that instance. A certified Gershgorin bound is approximately 1.13249 times the exact eigenvalue; using it online would reduce the step size and could increase the iteration count.

Because the present algorithm prioritizes computational efficiency, it retains finite power iteration and a 3\% empirical margin, while the exact-eigenvalue and Gershgorin computations remain offline audits. This choice follows the same principle as the one-pass node map: keep online work fixed and extremely small, state the approximation boundary precisely, and use task-level physical acceptance rather than an unsupported exact-optimality claim.

\section{Definitions of Residuals, Violations, and Timing Statistics}

\subsection{Outer Terminal Metrics}

To separate the stopping test from physical units, the terminal metric after nonlinear propagation is
\begin{equation}
e_f=\max\left\{
\frac{\norm{\vect r_N}_\infty}{\SI{500}{\meter}},
\frac{\norm{\vect v_N}_\infty}{\SI{50}{\meter\per\second}}
\right\}.
\end{equation}
It is evaluated on the nonlinearly propagated state, not on a scaled variable in the condensed equality. The relative control and time changes are
\begin{align}
e_q&=\frac{\norm{\vect q^{p+1}-\vect q^p}_2}
{\max(1,\norm{\vect q^p}_2)},\\
e_t&=\frac{|\tau^{p+1}-\tau^p|}
{\max(1,|\tau^p|)} .
\end{align}
They diagnose stabilization of the reference trajectory but do not replace $e_f$ or the physical feasibility tests. The virtual-buffer norm is reported separately so that a small condensed equality residual cannot conceal extensive borrowing of artificial feasibility.

\subsection{Inner Residuals}

Every fourth update evaluates the scaled equality residual
\begin{equation}
r_{\rm eq}=\norm{H\vect y-\vect b}_\infty
\end{equation}
and the relative fixed-point change
\begin{equation}
r_{\rm fp}=
\frac{\norm{D_c(\vect y^j-\vect y^{j-1})}_2}
{\max(1,\norm{D_c\vect y^j}_2)} .
\end{equation}
These quantities only describe stopping under the implemented map. Because that map is not the exact joint proximal operator, $r_{\rm fp}$ is not a complete KKT residual for the ideal composite problem. The offline audit separately evaluates primal feasibility, set violations, and the map error relative to converged Dykstra splitting.

\subsection{Task-Level Violations}

Nodewise validation records
\begin{align}
v_{\min}&=\max_k(T_{\min}-\norm{\vect T_k}_2)_+,\\
v_{\max}&=\max_k(\norm{\vect T_k}_2-T_{\max})_+,\\
v_{\theta}&=\max_k\left(
\cos\theta_{\max}\norm{\vect T_k}_2-
\vect e_y^\mathsf T\vect T_k\right)_+,\\
v_m&=\max_k(m_{\rm dry}-m_k)_+ .
\end{align}
It also checks mass monotonicity, agreement between the mass change and the trapezoidally integrated propellant consumption, and every nonlinear trapezoidal defect. The FOH high-order propagation reuses the same control nodes but applies a different integrator at dense query points, exposing intersample minimum-thrust violations that a nodewise test cannot detect.

\subsection{Latency Statistics}

P50, P95, and P99 are empirical quantiles of complete-plan samples; inner-call time, one-outer-step time, and batch throughput are not mixed into the same statistic. The 2000 nominal xPIPG repetitions describe the central timing distribution reliably. ECOS has only 200 repetitions, so its P99 is determined by approximately the two slowest observations and is treated as descriptive. For Monte Carlo results, repeated timings are first reduced to a within-case median and then aggregated across common successful cases. This prevents an operating-system interruption in a single repeat from being misinterpreted as a mission-envelope effect. Failed cases remain in the denominator of the success rate.

\section{Complete Algorithmic Sequence}

This appendix gives the complete sequence from initialization through timing and validation. The one-pass node operation is explicitly identified as an approximation. The word ``exact'' is reserved for block elimination of the linearized dynamics equalities.

\begin{algorithmblock}{Algorithm A1: Proximal-Regularized Condensed Sequential Convex Approximation}
\algline{Input:}{$N=31$, physical parameters, initial and terminal states, final-time bounds, $\lambda,w_q,w_t$, and stopping thresholds.}
\algline{Output:}{Nonlinear trajectory, control, final time, iteration history, and an independent validity flag.}
\algline{1.}{Linearly interpolate position, velocity, and mass; construct a control reference satisfying the pointing and magnitude bounds; set $\bar\tau=(15+45)/(2\times30)=1$.}
\algline{2.}{Evaluate the initial nonlinear metrics and set the terminal-buffer penalty to $\mu_\nu=50$.}
\algline{3.}{For outer index $p=1,\ldots,25$, perform Steps 4--13.}
\algline{4.}{Evaluate $\vect f,F_x,F_q$ at the 31 unique nodes and share node data between adjacent intervals.}
\algline{5.}{Use Algorithm A2 to assemble the trapezoidal interval models and condensed terminal sensitivities.}
\algline{6.}{Normalize the six terminal rows by \SI{500}{\meter} and \SI{50}{\meter\per\second}, apply unit-row-norm scaling, and append the $I/\sqrt{\mu_\nu}$ buffer columns.}
\algline{7.}{Set $\epsilon_{\rm in}=2\times10^{-3}$ initially and subsequently use $\max\{2\times10^{-4},\min(2\times10^{-3},0.05e_f)\}$.}
\algline{8.}{Construct the fuel-dominant strongly convex composite surrogate and call the latency-oriented xPIPG implementation in Algorithm A3.}
\algline{9.}{Record the control and time changes and accept the full update $\bar{\vect q}\leftarrow\vect q^\star,\bar\tau\leftarrow\tau^\star$.}
\algline{10.}{Using the preceding trajectory as the Newton initial guess, propagate exactly one full-step nonlinear implicit-trapezoidal trajectory.}
\algline{11.}{Recompute the objective, nonlinear terminal error, dry-mass margin, spectral estimate, step sizes, inner-iteration count, and buffer.}
\algline{12.}{Exit successfully when $e_f\le10^{-3}$ and the dry-mass check passes; retain the control and time changes as diagnostics.}
\algline{13.}{If $e_f>10^{-3}$, update $\mu_\nu\leftarrow\min(5\mu_\nu,2\times10^6)$.}
\algline{14.}{Apply Algorithm A4 to the final trajectory; never return an invalid candidate as a valid plan.}
\end{algorithmblock}

\begin{algorithmblock}{Algorithm A2: Trapezoidal Linearization, Schur Solve, and Active-Column Condensation}
\algline{Input:}{Reference trajectory $\bar{\vect x}_{1:N},\bar{\vect q}_{1:N},\bar\tau$.}
\algline{Output:}{$H_q,\vect h_\tau,\vect b$ and the interval linear models.}
\algline{1.}{Evaluate the baseline dynamics once per node; fill the position--velocity, density--altitude, mass, and control Jacobian entries analytically.}
\algline{2.}{Apply scaled forward differences only to the three velocity components, reusing the baseline dynamics.}
\algline{3.}{For every interval, form $M_k,R_k,U_k,V_k,\vect s_k$ and the reference trapezoidal defect.}
\algline{4.}{Form and factor the $3\times3$ velocity Schur complement of $M_k$ once; solve all right-hand sides in one batch, and reject the outer update if a pivot test fails.}
\algline{5.}{Recover $A_k,B_k,E_k,\vect g_k,\vect a_k$.}
\algline{6.}{Initialize the state offset, control sensitivity, and time sensitivity to zero.}
\algline{7.}{For $k=1,\ldots,N-1$, propagate the offset and time column while multiplying only the first $3(k+1)$ active control columns.}
\algline{8.}{Accumulate $B_k$ and $E_k$ in control blocks $k$ and $k+1$; retain future columns as structural zeros.}
\algline{9.}{Extract the six terminal position and velocity rows and combine them with the reference terminal error to obtain $H_q,\vect h_\tau,\vect b$.}
\algline{10.}{In debug validation, compare the condensed prediction with the complete linear dynamics chain for random $\delta\vect q,\delta\tau$.}
\end{algorithmblock}

\begin{algorithmblock}{Algorithm A3: Spectral-Step Latency-Oriented xPIPG Inner Loop}
\algline{Input:}{$P,\vect c,H,\vect b$, node reference directions, thrust thresholds, time bounds, and $\epsilon_{\rm in}$.}
\algline{Output:}{$\vect q^\star,\tau^\star,\vect\nu^\star$ and inner-loop diagnostics.}
\algline{1.}{Form $G=HH^\mathsf T$, perform 80 power iterations on the $6\times6$ matrix, and multiply by $1.0000001$. Treat the result as an empirical estimate, not a certified spectral upper bound.}
\algline{2.}{Set $L=\max_jP_{jj}$ and use $r=1,s=0.97$ to compute primal and dual steps $\alpha,\beta$.}
\algline{3.}{Cache $1-\alpha P_{jj}$, $-\alpha c_j$, and $\alpha H$; initialize the primal variables from the reference control and time.}
\algline{4.}{Choose six buffer variables analytically so that the initial condensed equality holds; apply column scales of 0.5 and 0.25 to the time and buffer columns, and set $\vect\xi=\vect y,\vect\eta=0$.}
\algline{5.}{For $j=1,\ldots,6000$, perform Steps 6--11.}
\algline{6.}{Compute the trial point $\vect u=\vect\xi-\alpha(P\vect\xi+\vect c+H^\mathsf T\vect\eta)$.}
\algline{7.}{At each control node, apply group shrinkage, cone--ball projection, and the minimum-thrust half-space projection in sequence. This single composition is the approximate node map.}
\algline{8.}{Clip time to $[15/30,45/30]$ and apply the smooth quadratic-gradient step to the virtual buffer.}
\algline{9.}{Evaluate $H\vect y$ once, update the temporary dual variable, and synchronously extrapolate primal and dual variables with $\rho_e=1.55$.}
\algline{10.}{Evaluate the equality residual and the scaled fixed-point change only on every fourth update or at the iteration limit; stop once both pass, provided at least 20 updates have elapsed.}
\algline{11.}{Update the cached $H\vect\xi$ through the extrapolation recurrence; do not copy the previous 100-vector on unchecked updates.}
\algline{12.}{Undo the column scaling, recover the physical buffer, and return the spectral estimate, step sizes, iteration count, and residuals.}
\end{algorithmblock}

\begin{algorithmblock}{Algorithm A4: Independent Physical Validation of a Candidate Trajectory}
\algline{Input:}{Problem data, final time, and the 31-node state and thrust sequences.}
\algline{Output:}{Validity flag and maximum violations by category.}
\algline{1.}{Check the solver status, finiteness of every array entry, and physical final-time bounds.}
\algline{2.}{Evaluate $\norm{\vect T_k}_2$ at every node and record lower- and upper-thrust violations.}
\algline{3.}{Evaluate $\cos\theta_{\max}\norm{\vect T_k}_2-\vect e_y^\mathsf T\vect T_k$ at every node to check the pointing cone.}
\algline{4.}{Check $m_k\ge m_{\rm dry}$ and nonincreasing mass at every node.}
\algline{5.}{Integrate $\norm{\vect T}_2/(I_{sp}g_0)$ with trapezoidal weights and compare it with $m_0-m_N$.}
\algline{6.}{Recompute every nonlinear trapezoidal defect using the final time to verify propagation consistency.}
\algline{7.}{Scale terminal position by \SI{500}{\meter} and velocity by \SI{50}{\meter\per\second}; require the six-dimensional infinity norm to be no greater than $10^{-3}$.}
\algline{8.}{Offline, propagate FOH thrust with high-accuracy \texttt{ode113} and report intersample violations; exclude this step from online timing.}
\algline{9.}{Return a valid online plan only when all task-level checks pass; printing, CSV generation, and plotting cannot alter the flag.}
\end{algorithmblock}

\begin{algorithmblock}{Algorithm A5: Complete-Plan Timing, Cross-Check, and Figure Generation}
\algline{1.}{Read the production C-PIPG and ECOS sources; write test copies, executables, and batch outputs only under \texttt{test}.}
\algline{2.}{Compile the benchmark executables with the archived Release options and record the floating-point mode and processor information.}
\algline{3.}{Warm up xPIPG for 20 end-to-end solves and ECOS for five; restart every solve from identical data and initialization.}
\algline{4.}{With printing and file output disabled, time initialization, outer approximation, inner solve, nonlinear propagation, and nodewise validation.}
\algline{5.}{Repeat the nominal xPIPG solve 2000 times and the ECOS solve 200 times; store iteration counts, terminal states, checksums, and quantiles.}
\algline{6.}{Read common Monte Carlo cases from the same CSV; count any Algorithm-A4 rejection as a failure without relaxing thresholds.}
\algline{7.}{For grid experiments, hold the physical and stopping parameters fixed and report complete latency, inner work, and physical outputs.}
\algline{8.}{Run offline proximal, exact-eigenvalue, interval-conditioning, and FOH continuous-time audits; exclude their time from the online result.}
\algline{9.}{Generate figures only from stored results and rerun the entire experiment whenever any frozen condition changes.}
\end{algorithmblock}

\section{Mapping from Paper Operations to Source Files}

The outer loop, $3\times3$ Schur condensation, the spectral estimate, the one-pass node map, and the xPIPG update are implemented in \path{landing_solver.c} in the C solver directory. Continuous dynamics, the semi-analytic Jacobian, implicit-trapezoidal propagation, and physical quantities are in \path{landing_dynamics.c}. The full-state reference uses \path{ecos_landing_solver.c} in the ECOS directory. Common entry points, benchmarks, shared cases, and figure postprocessing reside under \path{test}. The file \path{reviewer_audit.m} in the English-paper directory records the offline proximal, spectral, interval-conditioning, and continuous-time audits.

This mapping explains how the equations reach the implementation; it does not let source structure stand in for mathematical evidence. Block elimination establishes condensation equivalence. Offline audits quantify the proximal and spectral approximations. Task-level checks, independent of scaling and dual variables, determine whether a trajectory is accepted.